\providecommand{\href}[2]{#2}

\documentclass[12pt]{amsart}

\usepackage{graphicx}
\usepackage{amsthm}
\usepackage{amsmath, amssymb}
\usepackage{thmtools}
\usepackage{thm-restate}

\usepackage{hyperref}

\usepackage{cleveref}

\swapnumbers

\theoremstyle{plain}
\newtheorem{Thm}{Theorem}

\newtheorem{Lem}[Thm]{Lemma}

\theoremstyle{definition}
\newtheorem{Def}[Thm]{Definition}

\begin{document}
\begin{abstract}  

Haken showed that the Heegaard splittings of reducible 3-manifolds are reducible, that is, a reducing 2-sphere can be found which intersects the Heegaard surface in a single simple closed curve.    When the genus of the ``interesting'' surface increases from zero, more complicated phenomena occur.     Kobayashi showed that if a 3-manifold $M^3$ contains an essential torus $T$, then it contains one which can be isotoped to intersect a strongly irreducible Heegaard splitting surface $F$ in a collection of simple closed curves which are essential in $T$ and in $F$.  In general there is no global bound on the number of curves in this collection.    We give conditions under which a global bound can be obtained. 
\end{abstract}
\title{ Tori and Heegaard splittings}

\maketitle

 \author Abigail Thompson \footnote{Supported by the National Science Foundation grant 1207765 and  the Institute for Advanced Study, where the author was the  Neil Chriss and Natasha Herron Chriss Founders' Circle Member for 2015-2016.}

\section{Introduction}

In \cite{H}, Haken showed that the Heegaard splittings of reducible 3-manifolds are reducible, that is, a reducing 2-sphere can be found which intersects the Heegaard surface in a single simple closed curve.    When the genus of the ``interesting'' surface increases from zero, more complicated phenomena occur.  We explore conditions under which the picture remains simple when the manifold is irreducible but contains an essential torus.     

The motivation for this work is two-fold.    

First, Kobayashi \cite{K} showed that if a 3-manifold $M^3$ contains an essential torus $T$, then it contains one which can be isotoped to intersect a (strongly irreducible) Heegaard splitting surface $F$ in a collection of simple closed curves which are essential in $T$ and in $F$.  In general there is no global bound on the number of curves in this collection for an arbitrary genus Heegaard surface.    We give conditions under which a global bound exists.

Second, it is known (\cite{T}, \cite{He}) that if $M$ contains an essential torus $T$, then the {\it distance} of the Heegaard splitting, as defined by Hempel in \cite{He} is  at most 2.   So a Heegaard splitting of a toroidal manifold has distance at most 2, but of course a manifold with a distance 2 Heegaard splitting is not necessarily toroidal.   This naturally leads to the question:  given a distance 2 Heegaard splitting of a 3-manifold $M$, can we discern whether or not $M$ is toroidal?   We give a partial answer to this question.

Let $M^3$ be a  closed, orientable, irreducible 3-manifold.   Let $F$ be a minimal genus Heegaard surface for $M$, so $F$ splits $M$ into two handlebodies, $H_1$ and $H_2$.   Call a Seifert-fibered space with base space a disk and two exceptional fibers a {\it boundary small Seifert-fibered space}.

Our main theorem is:

\begin{Thm}\label{Main}
Let $M^3$ be a closed, orientable 3-manifold, and $F$ a minimal genus strongly irreducible  Heegaard splitting for $M$.    Let $T$ be an essential torus in $M$.   Then one of the following holds:\\

\noindent  i. There exists an essential surface $G$  and a minimal genus Heegaard surface $F'$ for $M$ such that $G$ intersects $F'$ in at most 4 essential simple closed curves.\\
ii. The minimal genus Heegaard decomposition of $M$ is not thin.   \\
iii. $M$ contains an essential torus bounding a boundary small Seifert-fibered space.\\
\end{Thm}

\subsection{Outline of the paper}

In Section 2 we give definitions and background information.  In Section 3 we prove Theorem \ref{Main}.

\section{Background and Definitions}

\subsection{Heegaard splittings and distance} Let $(H_1, H_2,F)$ be a Heegaard splitting of a closed orientable 3-manifold M, where $H_1$ and $H_2$ are handlebodies and $F = \partial{H_1} = \partial{H_2}$. The {\it genus} $g$ of the Heegaard splitting is the genus of the surface $F$. The Heegaard splitting is {\it reducible} if there exists an essential simple closed curve $c$ on $F$ such that $c$ bounds (imbedded) disks $D_1$ in $H_1$ and $D_2$ in $H_2$. The splitting is {\it stabilized} if there exist essential simple closed curves $c_1$ and $c_2$ on $F$ such that $c_i$ bounds an (imbedded) disk $D_i$ in $H_i$ and $c_1$ and $c_2$ intersect transversely in a single point. A stabilized splitting of genus at least 2 is reducible. The splitting is {\it weakly reducible} if there exist essential simple closed curves $c_1$ and $c_2$ on $F$ such that $c_i$ bounds an (imbedded) disk $D_i$ in $H_i$ for $i=1,2$ and $c_1$ and $c_2$ are disjoint. A splitting that is not weakly reducible is {\it strongly irreducible}. 

Hempel \cite{He} generalized the idea of strong irreducibility to define the {\it distance} of a Heegaard splitting to be the minimum length $r$ of a sequence $c_1,c_2,...,c_r$ of essential simple closed curves on $F$ such that $c_1$ bounds a disk in $H_1$, $c_r$ bounds a disk in $H_2$, and consecutive $c_i$'s are disjoint.    In this notation, a strongly irreducible Heegaard splitting has distance at least 1.   

\subsection{Thin position for 3-manifolds}  In \cite{S-T}, we define {\it thin position} for a closed, orientable 3-manifold $M$.     The idea is to replace a Heegaard splitting for $M$ with a different handle decomposition, which by some measure of complexity is potentially simpler than a Heegaard decomposition.   We include the basic definitions here; for more details, see \cite{S-T}.    

For $M$ a connected, closed, orientable 3-manifold, let $M=b_0\cup{N_1}\cup{T_1}\cup{N_2}\cup...\cup{N_k}\cup{T_k}\cup{b_3}$, where $b_0$ is a collection of 0-handles, $b_3$ is a collection of 3-handles, and for each $i=1,2,...,k$, $N_i$ is a collection of 1-handles and $T_i$ is a collection of 2-handles.   Let $S_i$, $i=1,2,...,k$, be the surface obtained from $\partial[b_0\cup{N_1}\cup{T_1}\cup...\cup{N_i}]$ by deleting all spheres bounding 0- or 3-handles in the decomposition.     The complexity of a (connected, closed, orientable) genus g surface $F$, $c(F)$,  is $2g-1$, and the complexity of a disconnected surface   is the sum of the complexities of its components.   Define the {\it width of the decomposition} of $M$ to be the set of integers $\{c(S_i)\}$, $i=1,2,...,k$.    We compare lists from two different decompositions using lexicographical ordering.     The {\it width} of $M$ is the the minimal width over all decompositions of $M$.     A handle decomposition for $M$ realizing the width of $M$ is called {\it thin position} for $M$.

It is straightforward to see that if a Heegaard splitting of $M$ is weakly reducible, then it is possible to re-arrange the handles of the splitting to obtain a thinner decomposition of $M$ than that provided by the Heegaard splitting.   What is less obvious is that it is possible for a minimal genus, strongly incompressible Heegaard splitting of $M$ to fail to be  thin position for the manifold.   This possibility arises in part (ii) of our main theorem.

\section{Proof of Theorem 1}

Let $M^3$ be a closed, orientable 3-manifold.   Let $F$ be a minimal genus strongly irreducible  Heegaard splitting for $M$.    Let $T$ be an essential torus in $M$.

By \cite{K}, we can isotop $T$ so that $T$ intersects $F$ in a collection $C$ of simple closed curves, each of which is essential both in $T$ and in $F$.   If the number of curves in $C$ is less than or equal to four, we are done, so assume that the number of curves in $C$ is at least six.   Hence $T$ is split by $C$ into at least six annuli.  We will use these annuli to obtain an annulus in $H_1$ which is disjoint from a ``good" pair of compressing disks in $H_2$.     

The proof of the theorem will follow from two lemmas.   The first produces a good pair  of compressing disks or an essential torus bounding a boundary small Seifert-fibered space.   The second uses a good pair of compressing disks to either produce the desired essential surface or to give a new, thinner,  handle decomposition of $M$.

\begin{Def} Let $H$ be a handlebody  and let $D$ and $E$ be disjoint compressing disks for $H$.  We say that $D$ and $E$  are {\em dependent} if either $D$ and $E$  are parallel in $H$,  or if $D_1$, say, cuts off a solid torus in which $D_2$ bounds a meridian disk.   We say the pair $(D,E)$ is {\em good} if at least one of $D$ and $E$ is non-separating and $D$ and $E$ are not dependent.       Suppose $(D,E)$ is a good pair of compressing disks for $H$, and let $F'$ be the boundary of the handlebody(ies) obtained by compressing $H$ along $D$ and $E$.   Note that $c(F')\leq c(F)-3$.  Indeed, the point of using a ``good" pair of disks is to ensure this drop in complexity.         

\end{Def}

\begin{Def}   Let $A$ be an annulus properly imbedded in a 3-manifold $M$.   Let $M'$ be obtained from $M$ by removing an open neighborhood of $A$.   We say $M'$ is obtained from $M$ by {\em surgering along $A$}.    In a slight abuse of notation, we also say that $\partial(M')$ is obtained from $\partial(M)$ by surgering along $A$.

\end{Def}

\begin{Lem}\label{LA}

Let $M^3$ be a closed, orientable 3-manifold.   Let $F$ be a minimal genus strongly irreducible  Heegaard splitting for $M$, splitting $M$ into the handlebodies $H_1$ and $H_2$.    Let $T$ be an essential torus in $M$.   Assume $T$ intersects $F$ in a collection $\bf{C}$ of simple closed curves which are essential on $T$ (and on $F$).   Assume the number of these curves has been minimized (among all choices of $T$ and $F$) and is greater than or equal to six. Then there exists a good pair of disks $(D,E)$ in $H_2$ disjoint from one of the annuli $A$ in $T\cap{H_1}$, or $M$ contains an essential torus $T$ bounding a boundary small Seifert-fibered space.

\end{Lem}

Proof of Lemma \ref{LA}:

Let $\bf{A}$ be the collection of annuli in $T\cap{H_2}$.   Every annulus in $\bf{A}$ is boundary compressible.    Find two annuli $B_1$ and $B_2$ in $\bf{A}$ so that $B_1$ can be boundary compressed disjoint from all other annuli in $\bf{A}$ to obtain the disk $D_1$  and then $B_2$ can be boundary compressed disjoint from all remaining annuli and $D_1$ to obtain the disk $D_2$.     Notice that $\partial{D_1}$ is disjoint from all curves in $\bf{C}$ and $\partial{D_2}$ is disjoint from all curves in $\bf{C}$ except possibly $\partial{B_1}$.   Since there are at least six curves in $\bf{C}$,  there is at least one annulus $A$ in $T\cap{H_1}$ with boundary disjoint from $\partial{D_1}\cup\partial{D_2}$.    

If the pair ${D_1,D_2}$ is good we are done.  

If ${D_1,D_2}$ is not good then either ${D_1}$ and ${D_2}$ are dependent or both ${D_1}$ and ${D_2}$ are separating, or possibly both.  

Suppose ${D_1}$ and ${D_2}$ are both separating but not parallel.   Then $\partial{D_1}\cup\partial{D_2}$  divides $F$ into three components.    At least one of these components is disjoint from $\partial{A}$.     Replace one of the $D_i$'s with a non-separating compressing disk $E$ for $H_2$  whose boundary is contained in the boundary of this component, making the selection to avoid dependence of the resulting pair.  Hence we can produce a good pair.     

Suppose ${D_1}$ and ${D_2}$ are dependent.  Then they are either parallel, or $D_1$, say, cuts off a solid torus in which $D_2$ bounds a meridian disk.    We can reconstruct the annuli $B_1$ and $B_2$ from $D_1$ and $D_2$ by attaching bands $d_1$ and $d_2$ to them.      

If ${D_1}$ and ${D_2}$ are parallel, then either  $B_1$ and $B_2$ are parallel in $H_2$ or $\partial{B_1}\cup\partial{B_2}$ bounds a 4-punctured sphere $P$ in $F$.      Note that $\partial{A}$ is disjoint from  $\partial{B_1}\cup\partial{B_2}$.   Some annulus (which by an abuse of notation we will call $A$) of $T\cap{H_1}$ must be boundary compressible in $H_1$ and the boundary compression must be disjoint from $\partial{B_1}\cup\partial{B_2}$.  Hence both boundary components of $A$ lie in $P$ or both lie outside of $P$.     If $\partial{A}$ lies completely outside $P$,  then boundary compressing $A$ yields a compressing disk for $H_1$ disjoint from $D_1$, contradicting strong irreducibility.     If $\partial{A}$ lies completely inside $P$, then at least two components of $\bf{C}$ must be parallel on $F$.  

Suppose $D_1$ cuts off a solid torus in which $D_2$ bounds a meridian disk.  Then reconstructing the annuli $B_1$ and $B_2$ from $D_1$ and $D_2$ by attaching bands $d_1$ and $d_2$ to them yields at least two of the curves in $\partial{B_1}\cup\partial{B_2}$ are parallel on $F$.    

In both cases, we obtain an annulus $S$ (the annulus of parallelism) on $F$ such that $\partial{S}$ lies on $T$.  Since $T$ was chosen to minimize the number of curves of intersection with $F$, $S$ is not parallel into $T$.   We can construct two new tori $T_1$ and $T_2$ by surgering $T$ along $S$, each of which can be isotoped to have fewer curves of intersection with $F$ than $T$.   Since $T$ was chosen to minimize the number of curves of intersection with $F$, both $T_1$ and $T_2$ must be inessential tori, hence (because $M$ is prime) each bounds a solid torus in $M$.   Then $T$ bounds a boundary small Seifert-fibered space.   
 
\begin{Lem}\label{LB}

Let $M^3$ be a closed, orientable 3-manifold.   Let $F$ be a minimal genus g strongly irreducible  Heegaard splitting for $M$, splitting $M$ into the handlebodies $H_1$ and $H_2$.   Let $A$ be an incompressible non-boundary-parallel annulus properly imbedded in $H_1$, and let $D$ and $E$ be a good pair of  compressing disks for $H_2$, such that $\partial{D}\cup\partial{E}$ is disjoint from $\partial{A}$.  Then at least one of the following holds: \\
1. There exists an essential surface $G$ that intersects $F$ in at most 4 essential simple closed curves.\\
2. The minimal genus Heegaard decomposition of $M$ is not thin.   \\
3. The surface $F'$ obtained by surgering $F$ along $A$ is also a Heegaard surface.

\end{Lem}

Proof of Lemma \ref{LB}:

Case 1:  $A$ is non-separating in $H_1$.

Let $H_1'$ be the manifold obtained from $H_1$ by surgering along $A$. Since $A$ is non-separating and incompressible, $H_1'$ is a handlebody of genus g.   Let $J$ be the complement of $H_1'$ in $M$.    If $J$ is a handlebody then possibility 3 holds and we are done.   

Assume $J$ is not a handlebody.     Since $\partial{D}\cup\partial{E}$ is disjoint from $\partial{A}$, $D$ and $E$ are compressing disks for $\partial{J}$.    Let $\bf{D}$ be a complete minimal collection of compressing disks for $J$ including $D$ and $E$ and let $L$ be the manifold obtained by compressing $J$ along $\bf{D}$.   Since $(D,E)$ is good, $c(\partial{L})\leq (2g-4)$.

Subcase A:  Some component $G$ of $\partial{L}$ is incompressible in $M$.     Then, by reconstructing $J$, we see that $G$ is an incompressible surface in $M$ and $G$ intersects $F$ in at most four essential simple closed curves.    

Subcase B:  Some component $G$ of $\partial{L}$ is compressible in $M$.  Since $\bf{D}$ is complete, $G$ is incompressible into $L$, hence it must be compressible into $M-L$.     By \cite{C-G}, the Heegaard splitting of $M-L$ given by $F'$ is weakly reducible, hence the width of $M-L$ is less than $2g-1$.   Starting with $\partial L$, however, we can complete the handle decompostion of $M$ by re-attaching $A$ and then completing the compressions from $H_2$.     Hence $L$ has width at most $2g-2$.     So the initial Heegaard splitting of $M$ was not thin position for $M$.     

Case 2:  $A$ is separating in $H_1$.

This case is similar, with a slightly more careful complexity count.   

Let $H_1'$ be the component obtained from $H_1$ by surgering along $A$ which contains $\partial{D}$ and $\partial{E}$ (since $F$ is weakly incompressible, both are in the boundary of one component, or else the boundary of one or the other would be disjoint from the disk obtained by boundary compressing $A$). $H_1'$ is a handlebody of genus at most $g$.   Let $J$ be the complement of $H_1'$ in $M$.    If $J$ is a handlebody then possibility 3 holds and we are done. 

If $J$ is not a handlebody, the argument proceeds as before.

\smallskip
\noindent {\em Conclusion of proof of Theorem \ref{Main}:}

Let $M^3$ be a closed, orientable 3-manifold.   Let $F$ be a minimal genus strongly irreducible  Heegaard splitting for $M$.    Let $T$ be an essential torus in $M$.   Assume $T$ intersects $F$ in a collection of simple closed curves which are essential on $T$ (and on $F$).   Among all such pairs $T$ and $F$, choose the pair with the fewest possible number of such curves of intersection.    If the number of curves is less than or equal to four, we are done, so assume the number of curves is at least six. Then by Lemma \ref{LA}, either we can find an essential torus bounding a boundary small Seifert-fibered space, or there exists a good pair of disks $(D,E)$ in $H_2$ disjoint from one of the annuli $A$ in $T\cap{H_1}$.   Applying Lemma \ref{LB}, either conditions (i) or (ii) hold, and we are done, or the surface $F'$ obtained by surgering $F$ along $A$ is also a Heegaard surface.   But then $F'$ is a minimal genus Heegaard surface intersecting $T$ in two fewer essential simple closed curves.   If $F'$ is a strongly irreducible Heegaard splitting, this contradicts our choice of $T$ and $F$.  If $F'$ is weakly reducible, then the Heegaard splitting given by $F'$ is not thin, hence the splitting given by $F$ is not thin.

 \

\begin{flushright}
Abigail Thompson\\ Department of Mathematics\\
University of California\\ Davis,
CA 95616\\ e-mail: thompson@math.ucdavis.edu\\
\end{flushright}

\end{document}